\documentclass[12pt,a4paper,reqno]{amsart}

\usepackage[utf8]{inputenc}
\usepackage{enumerate}
\usepackage[margin=2.5cm]{geometry}
\usepackage{graphicx} 
\usepackage{overpic}
\usepackage{hyperref}

\newtheorem{lema}{Lemma}[section]
\newtheorem{teo}[lema]{Theorem}
\newtheorem{pro}[lema]{Proposition}

\newtheorem{defini}[lema]{Definition}
\newtheorem{coro}[lema]{Corollary}

\DeclareMathOperator{\e}{e}

\newcommand{\R}{\ensuremath{\mathbb{R}}}

\title[Limit cycles in piecewise quadratic systems with an invariant straight line]{Bifurcation of limit cycles in piecewise quadratic differential systems with an invariant straight line}

\author[L. P. C. da Cruz] {Leonardo P. C. da Cruz}
\address{Departamento de Matem\'{a}tica, Universidade Federal de S\~{a}o Carlos, Rodovia Washington Lu{\'{\i}}s, Km 235, Caixa Postal 676, SI-2000, São Carlos (SP), Brazil}
\email{leonardo@mat.uab.cat}

\author[J. Torregrosa]{Joan Torregrosa}
\address{Departament de Matem\`{a}tiques, Universitat Aut\`{o}noma de Barcelona, 08193 Be\-lla\-ter\-ra, Barcelona (Spain); Centre de Recerca Matem\`atica, Campus de Bellaterra, 08193 Bellaterra, Barcelona (Spain)}
\email{torre@mat.uab.cat}

\subjclass[2010]{Primary 34C07, 34C23, 37C27}

\keywords{Center-focus, cyclicity, limit cycles, weak-focus order, Lyapunov quantities}

\begin{document}
	
\begin{abstract} 
We solve the center-focus problem in a class of piecewise quadratic polynomial differential systems with an invariant straight line. The separation curve is also a straight line which is not invariant. We provide families having at the origin a weak-foci of maximal order. In the continuous class, the cyclicity problem is also solved, being $3$ such maximal number. Moreover, for the discontinuous class but without sliding segment, we prove the existence of $7$ limit cycles of small amplitude.
\end{abstract}

\maketitle

\section{Introduction}

In past years, a big interest in the study of the dynamics of piecewise systems has emerged, due to the fact that many real phenomena can be modeled with this class of systems. For example, the existence and uniqueness of periodic orbits or the existence of a continuum of periodic orbits. These problems appear in many areas of research. In particular in electrical and mechanical engineering, in control theory, and even in the analysis of genetic networks. See for example \cite{AcaBonBro2011,BerBudCha2018}. 

\begin{figure}[h]
\begin{center}
\includegraphics{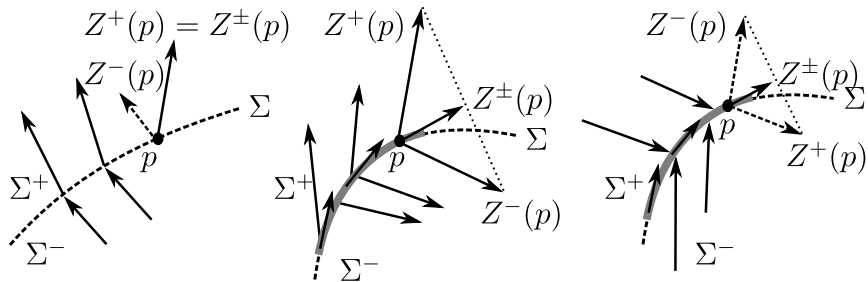}
\end{center}
\caption{Definition of the vector field on $\Sigma$ following Filippov's convention in the sewing, escaping, and sliding regions.}\label{fi:filipov}
\end{figure}

Usually, the simplest models are defined via planar piecewise polynomial vector fields $Z=(Z^+,Z^-)$ in the following way. Taking $0$ as a regular value of the function $h:\mathbb{R}^2\rightarrow \mathbb{R}$, we denote the discontinuity curve by $\Sigma=h^{-1}(0)$ and the two regions it delimits by $\Sigma^{\pm}=\{\pm h(x, y)>0 \}$. So, the piecewise vector field can be written as 
\begin{equation}\label{eq:1}
Z^{\pm}:(\dot{x},\dot{y})=(X^{\pm}(x,y),Y^{\pm}(x,y)),
 \text{ for } (x,y)\in \Sigma^{\pm},
\end{equation}
where $X^{\pm}$ and $Y^\pm$ are polynomials of degree $n$ in $\Sigma^\pm$.
The above piecewise vector field is continuous when it satisfies $Z^+=Z^-$ on the separation curve $\Sigma.$ Otherwise we will say that it is discontinuous. 
The local trajectories of $Z$ on $\Sigma$ was stated by Filippov in \cite{Fil1988} (see Figure~\ref{fi:filipov}). The points on $\Sigma$ where both vectors fields simultaneously point outward or inward from $\Sigma$ define the \emph{escaping} ($\Sigma^e$) and \emph{sliding region} ($\Sigma^s$), respectively. The interior of its complement on $\Sigma$ defines the \emph{crossing region} ($\Sigma^c$), and the boundary of these regions is constituted by tangential points of $Z^{\pm}$ with $\Sigma.$ As this work is restricted to the study of \emph{limit cycles of crossing type}, that we will refer to them only as \emph{limit cycles}, we do not recall here the precise definition of the vector field on $\Sigma^e$ and $\Sigma^s.$
Let $Z^{\pm}h$ denote the derivative of the function $h$ in the direction of the vector $Z^{\pm}$ that is, $Z^{\pm}h(p)=\langle \nabla h(p), Z^{\pm}(p)\rangle.$ Notice that $p\in\Sigma^c$ provided that $Z^+h(p)\cdot Z^-h(p) > 0,$ $p\in\Sigma^e\cup\Sigma^s$ provided that $Z^+h(p)\cdot Z^-h(p) < 0,$ and $p$ in $\Sigma$ is a tangential point of $Z^{\pm}$ provided that $Z^+h(p)Z^-h(p)=0.$ We say that $p \in \Sigma$ is a \emph{pseudo-equilibrium} of $Z,$ if $p$ is either a tangential point or an equilibrium of $Z^+$ or $Z^-.$ We call $p\in \Sigma$ an \emph{invisible fold} of $Z^+$ (resp. $Z^-$) if $p$ is a tangential point of $Z^+$ (resp. $Z^-$) and $(Z^+)^2h(p)<0$ (resp. $(Z^-)^2h(p)>0$). 

Let us consider that both differential systems in \eqref{eq:1} (when we thought them separately) have an equilibrium point at the origin such that the eigenvalues of their Jacobian matrices at $0$ have zero real part. For simplicity, we will consider only the cases when the linear part of each system is written in its normal form. Hence, after a time rescaling if necessary, \eqref{eq:1} writes as \begin{equation}\label{eq:3_1}
    Z^\pm=\begin{cases}
    \dot{x}=-y +\sum_{k=2}^{n} P^\pm_k(x,y),\\
    \dot{y}=x+\sum_{k=2}^{n} Q^\pm_k(x,y),
    \end{cases} \text{if} \ \ (x,y)\in \Sigma^{\pm} ,
\end{equation}
being $\Sigma^{\pm}=\{(x,y):\pm h(x,y)>0\}$ and $h$ a $\mathcal C^1$ function for which $0$ is a regular value. Obviously, for the problem that we would like to study, the time orientation is taken in order that the origin has a monodromic character. In this paper, we will assume that the discontinuity curve $\Sigma = h^{-1}(0)$ is a straight line passing through the origin. In fact, we will take them being one of the coordinates axes. As usual  $P^\pm_k,$ $Q^\pm_k$ denote homogeneous polynomials of degree $k.$ As in the analytic scenario, the problem of distinguishing whether the origin of \eqref{eq:3_1} is a center or a focus is also known as \emph{Poincar\'e center problem}, \emph{center-focus problem}, or just \emph{center problem}. In the piecewise polynomial class, there are other center problems that are not considered in this work. This is the case, for example, when we consider a pseudo-equilibrium point of fold-fold type. In addition to the difficulty of increasing the number of parameters, when we have fixed the degree of a vector field, we have to consider other types of centers appearing in the nonsmooth scenario. We will deal with this point in Section~\ref{se:suficientecondition}. A very related problem is the analysis of the number of limit cycles bifurcating from the origin.

It is well-known that linear vector fields have no limit cycles. But this is not the case in a piecewise scenario. Freire, Ponce, Rodrigo and Torres prove in \cite{FrePonRodTor1998} that only one limit cycle exists in continuous piecewise linear differential systems. For discontinuous piecewise linear differential systems Freire, Ponce, and Torres in \cite{FrePonTor2014} prove that only two limit cycles of small amplitude bifurcate from the origin. They also prove that a third (big) limit cycle exists. Also for this class of differential systems, the center problem from a monodromic equilibrium point but near the infinity is solved in \cite{FreiPonceTorre2021}, where the number of limit cycles bifurcating from the infinity is also considered. The authors prove that at least three limit cycles bifurcate from infinity. Both problems are related to using a transformation that moves the infinity to the origin. All the limit cycles are in fact \emph{crossing limit cycles}, because they cut the separation straight line $\Sigma$.

For quadratic vector fields, Bautin showed in \cite{Bautin1952} that the maximum number of limit cycles of small amplitude near an equilibrium point is three and, moreover, this upper bound is reached.
For quadratic discontinuous differential systems, this problem is studied in \cite{GouvTorre2021,JiLiuLi2022}. The work of Bautin is very appreciated because increasing the degree, these problems remain open. Thus, it is quite natural to restrict the study to some special families. For example, the planar quadratic vector fields have an invariant straight line. Cherkas, Zhilevich, and Rychkov, see \cite{CheZil1970, CheZil1972,Ryc1972}, proved that this family has only one limit cycle. For more details on this problem see \cite{CollLib1990} and \cite{Ye1986}. They prove that the canonical form of such systems is 
\begin{equation}\label{eq:0}
\begin{cases}
\dot{x}=-y+d x + l\,x^2+m\,x\,y+n\,y^2,\\
\dot{y}= x+b x\,y,
\end{cases}
\end{equation}
being $d\in [0,2),$ $m\geq0,$ and $b\ne0.$ We notice that, with these conditions, the origin is an equilibrium point of monodromic non-degenerate type. We will analyze this problem in the piecewise framework with 2 zones separated by a straight line passing through the origin. Hence, the dimension of the space of parameters $(d,l,m,n,b)\in \mathbb{R}^5$ will be doubled, $(d_1,d_2,l_1,l_2,m_1,m_2,n_1,n_2,b_1,b_2)\in \mathbb{R}^{10}.$ But, as we will see, the number of limit cycles will increase much more than doubled. A first natural extension is to consider that the invariant straight line is the same in both zones. Therefore we will assume $b_1=b_2=b.$ As in \cite{CollLib1990}, after the rescaling $(x,y)\rightarrow (- x/b,-y/b)$ if necessary, we can assume $b=-1.$ This change of variables does not modify the crossing limit cycles and the dimension of the parameter space decreases to $8$. 

An interesting phenomenon is that the number of limit cycles increases and it depends on the separation straight line. We will study two situations but fixing the canonical form \eqref{eq:0}: firstly when the separation straight line is the horizontal axis and secondly when it is the vertical axis. As we will see in the following results, the highest number of limit cycles is obtained for the second situation in the discontinuous case while in the continuous case is for the first one.

Consequently, the first piecewise quadratic differential system is
\begin{equation}\label{eq:1_1}
   Z_j=\begin{cases}
    \dot{x}=-y+d_j x + l_j\,x^2+m_j\,x\,y+n_j\,y^2,\\
    \dot{y}= x(1-y),
    \end{cases} \text{if} \ \ (x,y)\in \Sigma^{\mathcal{H}}_j,
\end{equation}
where the discontinuity straight line is $\Sigma^{\mathcal{H}}=\{(x,y): y=0\}$ and  $\Sigma^{\mathcal{H}}_j =\{(-1)^jy<0\},$ for $j=1,2.$ The second piecewise quadratic differential system is
\begin{equation}\label{eq:1_1a}
Z_j=\begin{cases}
		\dot{x}=-y+d_j x + l_j\,x^2+m_j\,x\,y+n_j\,y^2,\\
		\dot{y}= x(1-y),
	\end{cases} \text{if} \ \ (x,y)\in \Sigma^{\mathcal{V}}_j,
\end{equation}
where the discontinuity straight line is $\Sigma^{\mathcal{V}}=\{(x,y): x=0\},$ with $\Sigma^{\mathcal{V}}_j =\{(-1)^jx<0\},$ for $j=1,2.$ We notice that, in general, \eqref{eq:1_1} and \eqref{eq:1_1a} are discontinuous piecewise differential systems. The first main result (Theorem~\ref{thmmainn}) provides a lower bound for the number of limit cycles of small amplitude in both situations. It is remarkable that in the second one we are using all the parameters for having a complete unfolding. For the continuous cases, labeled as \eqref{eq:1_1c} and \eqref{eq:1_1ac} respectively, the local cyclicity problem is completely solved. See Theorems~\ref{Thm:continuo4} and \ref{Thm:continuo5}.

\begin{teo}\label{thmmainn} 
There are values of the parameters such that from the origin of systems \eqref{eq:1_1} and \eqref{eq:1_1a} bifurcate $4$ and $7$ crossing limit cycles of small amplitude, respectively, multiplicities taken into account.
\end{teo} 
As we will see in the proof of the above result, the limit cycles are obtained by studying which are the maximal orders of weak-foci together with the respective unfoldings. Both numbers provide lower bounds for the local cyclicity of the origin. For the study of the corresponding upper bounds, the cyclicity of the centers should be also studied. Which needs a more accurate analysis. This problem is completely solved in the next two results, getting the least upper bounds for the maximum number of crossing limit cycles of small amplitude that can bifurcate from the origin, multiplicities taken into account, that is, providing its \emph{local cyclicity}. We remark the large difference between the number of limit cycles taking into account that both families have the same number of parameters. We also observe that the considered families have no sliding segment near the origin. As we will explain in Section~\ref{se:preliminaries}, the pseudo-Hopf bifurcation does not take place. If we were interested in this kind of bifurcation, we could add an extra parameter in the first components of system~\eqref{eq:1_1a} obtaining an extra limit cycle. We notice that the study of this phenomenon in system~\eqref{eq:1_1} breaks the chosen canonical form and also the invariant straight lines.

\begin{teo}\label{Thm:continuo4}
The piecewise differential equation~\eqref{eq:1_1} is continuous if and only if $l_1=l_2=:l$ and $d_1=d_2=:d.$ Therefore, it becomes 
\begin{equation}\tag{$\ref{eq:1_1}^c$}\label{eq:1_1c}
\begin{cases}
\dot{x}=-y+d \,x + l \,x^2+m_j\,x\,y+n_j\,y^2,\\
\dot{y}= x(1-y),
\end{cases} \text{if} \ \ (x,y)\in \Sigma^{\mathcal{H}}_j,
\end{equation} 	
being $\Sigma^{\mathcal{H}}_j =\{(-1)^jy<0\},$ for $j=1,2.$ Moreover, the local cyclicity of the origin of the above family is $3$, multiplicities taken into account.
\end{teo}

\begin{teo}\label{Thm:continuo5}
The piecewise differential equation~\eqref{eq:1_1a} is continuous if and only if $n_1=n_2=:n.$ Therefore, it writes as
\begin{equation}\tag{$\ref{eq:1_1a}^c$}\label{eq:1_1ac}
\begin{cases}
\dot{x}=-y+d_j\,x + l_j\,x^2+m_j\,x\,y+n\,y^2,\\
\dot{y}= x(1-y),
\end{cases} \text{if} \ \ (x,y)\in \Sigma^{\mathcal{V}}_j,
\end{equation}	
being $\Sigma^{\mathcal{V}}_j =\{(-1)^jx<0\},$ for $j=1,2.$ Moreover, the local cyclicity of the origin of the above family is $2$, multiplicities taken into account.	
\end{teo}

After a detailed analysis of the number of limit cycles of small amplitude bifurcating from the origin in families \eqref{eq:1_1} and \eqref{eq:1_1a}, we finish providing an answer to the respective center-focus problems, firstly for the discontinuous case and secondly for the continuous one.
 
\begin{teo}\label{centro1} For family \eqref{eq:1_1} with $d_1=d_2=0,$ the origin is a center if, and only if, one of the next conditions holds:
\begin{itemize}
	\item[$(\mathcal{H}_1)$] $m_1=m_2= 0;\vspace{0.1cm}$
	\item[$(\mathcal{H}_2)$] $l_1-l_2=m_1-m_2=l_1+n_2=l_1+n_1=0;\vspace{0.1cm}$	\item[$(\mathcal{H}_3)$] $l_1+l_2+1=m_1-m_2= n_1 +n_2-1=0.\vspace{0.1cm}$
\end{itemize}	
\end{teo}

\begin{teo}\label{centro2} For family \eqref{eq:1_1a} with $d_1=d_2=0,$ the origin is a center if, and only if, one of the next conditions holds:
\begin{itemize}
	\item[$(\mathcal{V}_1)$] $l_1-l_2=m_2+m_1=n_1-n_2=0;\vspace{0.1cm}$
	\item[$(\mathcal{V}_2)$] $l_1-l_2=n_2+l_2=n_1+l_2=0;\vspace{0.1cm}$
    \item[$(\mathcal{V}_3)$] $l_1 -2l_2=m_2=n_2=n_1 +2l_2=0;\vspace{0.1cm}$
    \item[$(\mathcal{V}_4)$] $2l_1-2l_2+1=m_2=m_1=n_2= n_1-1=0;\vspace{0.1cm}$
    \item[$(\mathcal{V}_5)$] $2l_1-2l_2-1=m_2=m_1=n_1=n_2-1=0;\vspace{0.1cm}$
    \item[$(\mathcal{V}_6)$] $l_1-2l_2-1=m_2=n_2-1=n_1+2l_2+1=0;\vspace{0.1cm}$
    \item[$(\mathcal{V}_7)$] $2l_1-l_2=m_1=n_1=n_2+l_2=0;\vspace{0.1cm}$
    \item[$(\mathcal{V}_8)$] $2l_1-l_2+1=m_1=n_2+l_2=n_1-1=0.$
\end{itemize}	
\end{teo}

The next corollaries follow straightforwardly from the above results, 
considering the continuity conditions given in Theorem~\ref{Thm:continuo4} when the corresponding family has an equilibrium point of weak-focus type at the origin.
\begin{coro}\label{coroh} For family \eqref{eq:1_1c} with $d=0,$ the origin is a center if, and only if, one of the next conditions holds:
	\begin{itemize}
		\item[$(\mathcal{H}_1^c)$] $m_1=m_2= 0;\vspace{0.1cm}$
		\item[$(\mathcal{H}_2^c)$] $m_1-m_2=l+n_2=l+n_1=0;\vspace{0.1cm}$
		\item[$(\mathcal{H}_3^c)$] $2l+1=m_1-m_2= n_1 +n_2-1=0.$	
	\end{itemize}	
\end{coro}

\begin{coro}\label{centvc} For family \eqref{eq:1_1ac} with $d_1=d_2=0,$ the origin is a center if, and only if, one of the next conditions holds:
	\begin{itemize}
		\item[$(\mathcal{V}_1^c)$] $l_1-l_2=m_2+m_1=0;\vspace{0.1cm}$
		\item[$(\mathcal{V}_2^c)$] $l_1-l_2=n+l_2=n+l_2=0;\vspace{0.1cm}$
		\item[$(\mathcal{V}_3^c)$] $l_1=l_2=m_2=n=0;\vspace{0.1cm}$
		\item[$(\mathcal{V}_6^c)$] $l_1+1=l_2+1=m_2=n-1=0;\vspace{0.1cm}$
		\item[$(\mathcal{V}_7^c)$] $l_1=l_2=m_1=n=0;\vspace{0.1cm}$
		\item[$(\mathcal{V}_8^c)$] $l_1+1=l_2+1=m_1=n-1=0.$
	\end{itemize}	
	
\end{coro}

The paper is structured as follows. In Section~\ref{se:preliminaries}, we present the basic tools necessary to prove the results of this work. In Section~\ref{se:suficientecondition}, we provide sufficient conditions so that the presented families have a center at the origin. Next, in Section~\ref{se:analytic}, we analyze the highest-order weak-foci equilibrium points together with the small amplitude limit cycle bifurcation and, as usual, the necessary conditions to have a center equilibrium point. The last section is devoted to showing which is the upper bound for the cyclicity when the families are continuous, finishing with the proofs of Theorems~\ref{Thm:continuo4} and \ref{Thm:continuo5}.

\section{The degenerate Hopf bifurcation}\label{se:preliminaries}
The main results of this paper follow studying the return map near an equilibrium point of monodromic type located in the separation straight line. In fact, studying the composition of two half-return maps because the proofs are mainly based on considering two piecewise polynomial vector fields having a nondegenerate equilibrium point of center-focus type. Hence, the analysis of both maps can be realized by computing the Taylor series of the solution with respect to the initial condition, working in polar coordinates. But, instead of using the composition of both maps, we will compute the difference map that is equivalent. We recall first how these Taylor series can be computed and then how they are used to study lower and upper bounds for the cyclicity. That is, the number of limit cycles of small amplitude bifurcating from the equilibrium point.

\bigskip

\begin{figure}[h]
\includegraphics[height=3.5cm]{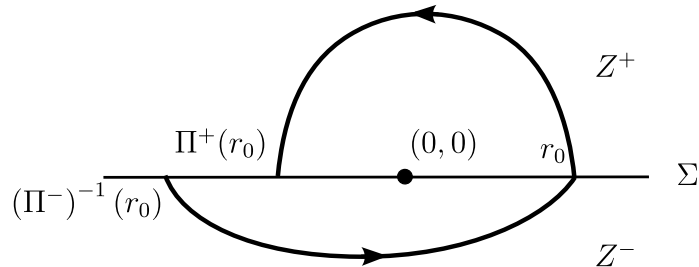}
\caption{The positive and negative half-return maps $\Pi^+$ and $\left(\Pi^-\right)^{-1},$  respectively.  }\label{fi:retunmap}
\end{figure}

The piecewise system \eqref{eq:3_1} can be written in polar coordinates, $(x,y)=(r\cos\theta,r\sin\theta),$ as 
\begin{equation}\label{eq:4}
 \frac{\mathrm{d} r^\pm}{\mathrm{d} \theta}\!=\!\frac{\sum_{k=2}^{n}R^\pm_k(\theta)r^k }{1+\sum_{k=2}^{n}\Theta^\pm_k(\theta)r^{k-1}}=\sum_{k=2}^{\infty}S^\pm_k(\theta)r^k, 
\end{equation} 
where 
\begin{equation*}
\begin{aligned}
R^\pm_k(\theta)=&\cos\theta \, P^\pm_k(\cos\theta,\sin\theta)+\sin\theta\,  Q^\pm_k(\cos\theta,\sin\theta),\\
\Theta^\pm_k(\theta)=&\cos\theta \, Q^\pm_k(\cos\theta,\sin\theta)-\sin\theta\,  P^\pm_k(\cos\theta,\sin\theta),
\end{aligned}
\end{equation*}
being $P^\pm_k,Q^\pm_k,$ and $S^\pm_k$ polynomials in $\sin\theta$ and $\cos\theta.$ We consider the solution of the initial value problems defined by \eqref{eq:4} with $r^+(r_0,0)=r_0$ and $r^-(r_0,\pi)=r_0,$ written in Taylor series with respect to $r_0,$ defined when $|r_0|\ll 1,$ as 
\begin{equation*}
    r^\pm(\theta,r_0)=\begin{cases}
    r_0+\sum_{k=2}^{\infty}u^+_k(\theta)r_0^k,\ \ \text{if} \ \ \theta\in (0,\pi), \\
    r_0+\sum_{k=2}^{\infty}u^-_k(\theta)r_0^k, \ \ \text{if} \ \ \theta\in (\pi,2\pi).
    \end{cases}   
\end{equation*}
Hence, we can define the positive half-return map $\Pi^+(r_0)=r^+(r_0,\pi)$ and the negative half-return map $\Pi^-(r_0)=r^-(r_0,2\pi).$ Instead of considering the composition of both maps we will define its equivalent displacement map, see for example \cite{CollPhoGassu1999},
\begin{equation}\label{eq:11}
\Delta(r_0)=\left(\Pi^-\right)^{-1}(r_0)-\Pi^+(r_0)=\sum_{k=2}^{\infty}W_kr_0^k,
\end{equation}
where, the function $\left(\Pi^-\right)^{-1}(r_0)$ is the inverse of the negative half-return map $\left(\Pi^-\right)(r_0),$ as it is illustrated in Figure~\ref{fi:retunmap}. As usual in this kind of analysis, the first nonvanishing $W_k$ is called the $k$th-order Lyapunov quantity of the piecewise polynomial system \eqref{eq:3_1}. This approach was also used in \cite{GassuTorre2003,HanZha2010,LiaHan2012}. 
It is well-known that, as the usual Lyapunov quantities for analytic vector fields, the coefficients $W_k$ are polynomials depending on the coefficients of the polynomial perturbations $P_k^\pm,Q_k^\pm.$ Finally, we will say that the origin of \eqref{eq:3_1} is a weak-focus of order $\ell$ when $W_j=0,$ $1\leq j\leq \ell-1$ and  $W_\ell\neq0.$ Moreover, the stability of the equilibrium point is given by the sign of the first nonzero Lyapunov quantity $W_\ell.$ 

This first non-vanishing coefficient provides also the stability of the equilibrium. As usual, choosing adequately the perturbation parameters in $P_k^\pm,Q_k^\pm$ we obtain limit cycles of small amplitude bifurcating, in this case, from the origin. This perturbation analysis is known as the degenerated Hopf bifurcation. 

The classical Hopf bifurcation in the study of analytic planar vector fields is characterized by the birth of a limit cycle from a weak-focus of first-order. More specifically, in the analytical context, the first nonvanishing coefficient has always an odd subscript and the limit cycle, which is of small amplitude, bifurcates from the origin changing adequately the sign of the trace of the Jacobian matrix of the corresponding system near the equilibrium point. See more details in \cite{AndLeoGorMai1973,Roussarie1998}.

This bifurcation, in the piecewise context, is associated with the study of the return map in a fold-fold type point. In this case, the first non-vanishing coefficient has always an even subscript. The limit cycle of small amplitude bifurcates from the origin changing also the stability of the origin. In this case, the size of the sliding segment takes the role of the trace in the analytic context. See more details of this phenomenon of codimension $1$ in \cite{Kuznetsov2003} or also in \cite{GuarSeaTei2011} for the codimension $2$ case. Recently, this bifurcation is also known as pseudo-Hopf bifurcation, see \cite{CasLliVer2017}.

When, as in our case, we fix the existence of an equilibrium of monodromic type in both upper and lower systems of \eqref{eq:3_1}, all the coefficients in $r_0$ appear in the Taylor development \eqref{eq:11}. All our families \eqref{eq:1_1}, \eqref{eq:1_1c}, \eqref{eq:1_1a}, and \eqref{eq:1_1ac} have the origin as an equilibrium point. Hence, $W_0=0$ and we have no sliding (nor escaping) segment. 

It is easy to check that by adding the trace parameter in upper and lower systems \eqref{eq:3_1} generically the first non-vanishing coefficient in \eqref{eq:11} is $W_1.$ As in the analytic scenario when we have an equilibrium point of focus type. In fact, $W_1=0$ if and only if the sum of the traces (of the upper and lower systems) is zero. Clearly, when $W_1=0$ and $W_2\neq0,$ the stability of the equilibrium point is given by the sign of $W_2.$ More concretely, the origin of the system is stable (resp. unstable) when $W_2<0$ (resp. $W_2>0$). Consequently, in the system, when $W_1$ is a small enough positive (resp. negative) real number, a small (resp. unstable) limit cycle bifurcates from the origin. In this case, it is important to remark that the equilibrium point (or equilibrium points) remains located at the origin. So, this bifurcation is also similar to the one previously denoted as Hopf bifurcation.

From the above analysis, in the bifurcation of an analytic planar piecewise vector field, when we have a weak-focus of order $k$ we get (generically) $k$ limit cycles. See more details in \cite{GouvTorre2021}. This bifurcation problem with varying parameters and taking into account multiplicities is studied in \cite{HanYan2021}. As we will see, as our families are polynomial, the study of the complete unfolding is more intricate. We notice again that, as all our families \eqref{eq:1_1}, \eqref{eq:1_1c}, \eqref{eq:1_1a}, and \eqref{eq:1_1ac} have no sliding, we will only get up to $k-1$ limit cycles of small amplitude. Although for families \eqref{eq:1_1ac} and \eqref{eq:1_1a} this pseudo-Hopf bifurcation makes sense. Because the invariant straight line $1-y=0$ remains unchanged when a constant term is added in the first components.
	
\section{The sufficient conditions for the center problem}\label{se:suficientecondition}	
	
This section is devoted to proving that the families in Theorems~\ref{centro1} and \ref{centro2} are centers. This is done in Propositions~\ref{cent3} and \ref{cent4}, respectively. We consider centers such that the period annulus is formed only by crossing periodic orbits. The key point is based, except by a special change of variables introduced in Definition~\ref{twin}, on the existence of three centers type. They are rigorously defined in Theorem~\ref{teocentros}. The first ones are of Darboux type because they have a piecewise first integral; the second ones have the usual time-reversibility, and the third ones are also symmetric but have the identity as the half-return map. 
	
To simplify the reading, in the next definition and the main result we take $\Sigma^{\mathcal{H}}=\{y=0\}$ as the separation curve. The result and the definition can be easily generalized considering other separation curves.
 
\begin{defini}\label{twin} Let $\Psi^\pm:\R^2\rightarrow\R^2$ be bijective transformations and $\Sigma^{\mathcal{H}}=\{y=0\}$. We say that \[
	\Psi(x,y)=\begin{cases}\Psi^+(x,y),& \text{ if }y>0,\\
	\Psi^-(x,y),& \text{ if }y<0,\end{cases}
	\] is a twin $\Sigma^{\mathcal{H}}$-transformation when $\Psi^+(x,0)=\Psi^-(x,0).$
\end{defini}

\begin{figure}[h]
		\includegraphics[height=3.5cm]{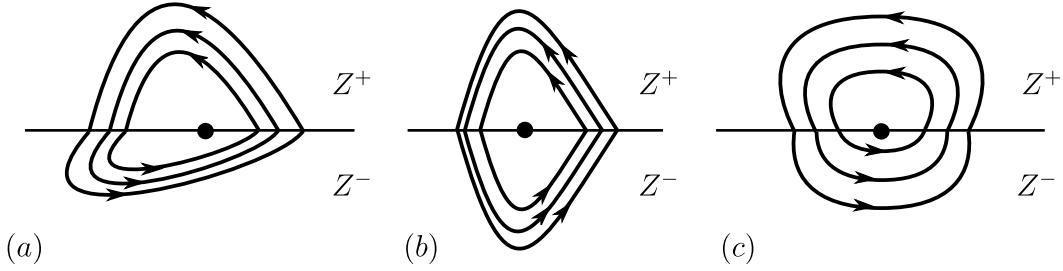}
	\caption{The three different center types detailed in Theorem~\ref{teocentros}}\label{fi:centros}
\end{figure}

\begin{teo}\label{teocentros} Let $Z$ be a piecewise differential system of the form \eqref{eq:3_1} with $\Sigma^{\mathcal{H}}=\{y=0\}.$ Then, applying a twin $\Sigma^{\mathcal{H}}$-transformation if necessary, we have a center at the origin in the following cases: 
\begin{enumerate}[(a)] 
 	\item There exist first integrals $H^\pm$ of $Z^\pm$ satisfying $H^+(x,0)=H^-(x,0).$   	
 	\item $Z$ is invariant with respect to the change
 	\begin{equation}\label{m1}
 	(x,y,t)\rightarrow (x,-y,-t).
 	\end{equation} 
    \item $Z^\pm$ are invariant with respect to the change
     \begin{equation}\label{m2}
     (x,y,t)\rightarrow (-x,y,-t).
     \end{equation}
     See all the different cases drawn in Figure~\ref{fi:centros}.
 \end{enumerate}   
\end{teo}

We notice that for the third class of systems in the above result also the complete vector field $Z$ has the detailed symmetry.

\begin{proof}[Proof of Theorem~\ref{teocentros}]
We notice that from the piecewise differential system \eqref{eq:3_1}, it is easy to see that the origin is a monodromic point and this property will be used along the proof.

\begin{enumerate}[(a)] 
\item From the canonical form, as it is written the piecewise differential system \eqref{eq:3_1}, it is easy to check that the Taylor series of the first integrals should start as $H^\pm(x,y)=x^2+y^2+\cdots.$ Let $x_2<0<x_1$ small enough real numbers such that $(x_j,0)\in\Sigma^{\mathcal{H}},$ $j=1,2,$ and are connected by monodromic solutions of \eqref{eq:3_1} in upper and lower half-planes. Consider the function
\begin{equation*}
\widehat{H}(x_1,x_2)=\dfrac{H^\pm(x_2,0)-H^\pm(x_1,0)}{x_2-x_1}=\sum_{i=1}^{n}\widehat{H}_i(x_1,x_2),  
\end{equation*}    
being $\widehat{H}_i$ homogeneous polynomials of degree $i$ in $(x_1,x_2).$ We notice that we have only one function $\widehat{H}$ because, by hypothesis, $H^\pm(x,0)=H^\pm(x,0).$ Moreover, $\widehat{H}_1(x_1,x_2)=x_1+x_2$ and, when $(x_1,0)$ and $(x_2,0)$ are on the same level curve of $H^+,$ or $H^-$, we have $\widehat{H}(x_1,x_2)=0.$ As, near the origin $\partial\widehat{H}/\partial x_2\neq0,$ we can apply the Implicit Function Theorem to show the existence of a unique half-return map $g:I\rightarrow I,$ where $I$ is an interval containing the origin. Moreover, $g$ satisfies $\widehat{H}(x,g(x))=0$ and $g(x)=-x+\cdots.$ The proof follows from the uniqueness of the half-return map, because it is the same in both,  upper and lower, regions. See Figure~\ref{fi:centros}.$(a).$  

\item The monodromy property together with the time-symmetry \eqref{m1}, allows us to use the classical result of analytic reversibility systems to prove this item. See more details in \cite{LamRob1998}. A drawing of this situation can be seen in Figure~\ref{fi:centros}.$(b).$

\item The proof follows by the reversibility property \eqref{m2} that satisfy the upper and lower solutions. Both are symmetric with respect to the $y$-axis, so the respective half-return maps are equal. More concretely, each point $(x,0)$ in $\Sigma^{\mathcal{H}}$ is sent to the corresponding symmetric one $(-x,0).$ See this property in Figure~\ref{fi:centros}.$(c).$
\end{enumerate}
\end{proof}

The following results are direct consequences of the last theorem.

\begin{pro}\label{cent3} For each family $\mathcal{H}_i,$ $i=1,2,3$, listed in Theorem~\ref{centro1}, the corresponding piecewise system \eqref{eq:1_1} has a center at the origin.       
\end{pro}
\begin{proof} 
The first family $\mathcal{H}_1$ is a center using directly Theorem~\ref{teocentros}.$(c)$. The second family $\mathcal{H}_2$ satisfies that $Z_1=Z_2.$ Hence it is, in fact, a quadratic vector field of Lotka--Volterra type following the classification in \cite{Zol1994}. So there exists a first integral that coincides in both regions $\Sigma^{\mathcal{H}}_1$ and $\Sigma^{\mathcal{H}}_2.$ The proof finishes applying Theorem~\ref{teocentros}.$(a)$. The proof that the last family $\mathcal{H}_3$ has a center follows from Theorem~\ref{teocentros}.$(b),$ proving that, after writing the systems in the usual polar coordinates $(x,y)=\left(r\cos\theta,r\sin\theta\right),$ the change $r=4R/(4+A_j(\theta)R),$ with $j=1,2$ where 
\begin{equation*}
\begin{aligned}
A_1(\theta)=&(3-l_1-3n_2)\sin\theta+m_2\cos\theta-(1+l_1-n_2)\sin 3\theta-m_2\cos 3\theta,\\
A_2(\theta)=&(1+l_1+3n_2)\sin\theta+m_2\cos\theta+(1+l_1-n_2)\sin 3\theta-m_2\cos 3\theta
\end{aligned}
\end{equation*}
is a twin $\Sigma^{\mathcal{H}}$-transformation.
\end{proof}
We notice that both rational changes of coordinates of the above proof are known because they allow us to change a system with a homogeneous nonlinearity to an Abel differential equation, see for example \cite{Che1976}.

\begin{pro}\label{cent4} 
 For each family $\mathcal{V}_i,$ $i=1,\ldots,8$, listed in Theorem~\ref{centro2}, the corresponding piecewise system \eqref{eq:1_1a} has a center at the origin.       	
\end{pro}

\begin{proof} As the separation line is $\Sigma^{\mathcal{V}}=\{x=0\},$ we can apply Theorem~\ref{teocentros} after changing the variables $(x,y)$ by $(y,x).$	
	
The first family $\mathcal{V}_1$ is time-reversible with respect to the change  $(x,y,t)\rightarrow (-x,y,-t)$ and, from the comment above, we have a center applying Theorem~\ref{teocentros}.$(b)$. The proof for the remaining families follows using Theorem~\ref{teocentros}.$(a)$ and all the first integrals $H^\pm,$ here denoted by $H_1,H_2$, will be of Darboux type and they will write as
\[
H_j(x,y)= (f_{j,1}(x,y))^{\lambda_{j,1}}(f_{j,2}(x,y))^{\lambda_{j,2}}(f_{j,3}(x,y))^{\lambda_{j,3}}.
\]
Where $\lambda_{j,3}$ can be zero when only two invariant algebraic curves are necessary in the center characterization. It is necessary to check that all of them are well defined in a neighborhood of the origin. We only provide the polynomials $f_{j,1},f_{j,2},f_{j,3}$ and the exponents $\lambda_{j,1},\lambda_{j,2}, \lambda_{j,3}.$ Using the first remark of the proof, in addition of finding the first integrals we will check the condition  $H_1(0,y)=H_2(0,y),$ for all $y$, being $H_1$ and $H_2$ the first integrals defined in $x>0$ and $x<0$, respectively.

$\bullet$ For the case $\mathcal{V}_2$ we have $H_j(0,y)=(1-y)^{l_2}(l_2y+1),$ for $j=1,2.$ Being
\begin{equation*}
\begin{aligned}
f_{j,1}& =1-y,\\
f_{j,2}& =(-m_j+(4l_2^{2}+m_j^{2}+4l_2)^{1/2}) x/2+l_2y+1,\\
f_{j,3}& =(-m_j-(4l_2^{2}+m_j^{2}+4l_2)^{1/2}) x/2+l_2y+1,\\
\lambda_{j,1}& = l_2,\\
\lambda_{j,2}&=(1+m_j(4l_2^2+m_j^2+4l_2)^{-1/2})/2,\\
\lambda_{j,3}&=(1-m_j(4l_2^2+m_j^2+4l_2)^{-1/2})/2.
\end{aligned}
\end{equation*}

$\bullet$ For family $\mathcal{V}_3$ we can take
\begin{equation*}
\begin{aligned}
f_{1,1}& =f_{2,1}=1-y,\\
f_{1,2}& = (-m_1+(16l_2^{2}+m_1^{2}+8l_2)^{1/2}) x/2+2l_2y+1,\\
f_{1,3}& = (-m_1-(16l_2^{2}+m_1^{2}+8l_2)^{1/2}) x/2+2l_2y+1,\\
\lambda_{1,1}& =\lambda_{2,1}= 2l_2,\\
\lambda_{1,2}&=(1+2m_1(16l_2^2+m_1^2+8l_2)^{-1/2})/2,\\
\lambda_{1,3}&=(1-2m_1(16l_2^2+m_1^2+8l_2)^{-1/2})/2,\\
f_{2,2}& =-(2l_2+1)l_2x^2+2l_2y+1,\\
\lambda_{2,2}&=1.
\end{aligned}
\end{equation*}
Here $H_j(0,y)=(1-y)^{2l_2}(2l_2y+1),$ for $j=1,2.$ 

$\bullet$ For family $\mathcal{V}_4,$ we have $H_j(0,y)=(1-y)^{2l_2}(2l_2y+1),$ for $j=1,2,$  where
\begin{equation*}
\begin{aligned}
f_{1,1}&=-f_{2,1} =y-1,\\
f_{1,2}& =(2l_2+1)l_2x^2+(2l_2y+1)(y-1),\\
\lambda_{1,1}&=\lambda_{2,1}-1 =2l_2-1,\\
\lambda_{1,2}&=\lambda_{2,2} =1,\\
f_{2,2}& = -(2l_2+1)l_2x^2+2l_2y+1.\\
\end{aligned}
\end{equation*}

$\bullet$ For case $\mathcal{V}_5$ we have 
\begin{equation*}
\begin{aligned}
f_{1,1}&=f_{2,1} =y-1,\\
f_{1,2}& =(2l_2^2+3l_2+1)x^2-(2l_2+1)y-1,\\
\lambda_{1,1}&=\lambda_{2,1}+1 =2l_2+1,\\
\lambda_{1,2}&=\lambda_{2,2}=1,\\
f_{2,2}& = (2l_2^2+3l_2+1)x^2+((2l_2+1)y+1)(y-1),\\
\end{aligned}
\end{equation*}
with $H_j(0,y)=\left((2l_2+1)y+1\right)(y-1)^{2l_2+1},$ for $j=1,2.$

$\bullet$ In family $\mathcal{V}_6$ we have $H_j(0,y)=\left((2l_2+1)y+1\right)(y-1)^{2l_2+1},$ with $j=1,2,$ and 
\begin{equation*}
\begin{aligned}
f_{1,1}&=f_{2,1} =y-1,\\
f_{1,2}& =\big(-m_1+(16l_2^2+m_1^2+24l_2+8)^{1/2}\big)x/2+(2l_2+1)y+1,\\
f_{1,3}& =\big(-m_1-(16l_2^2+m_1^2+24l_2+8)^{1/2}\big)x/2+(2l_2+1)y+1,\\
\lambda_{1,1}&=\lambda_{2,1}+1 =2l_2+1,\\
\lambda_{1,2}&=[1+m_1(16l_2^2+m_1^2+24l_2+8)^{-1/2}]/2,\\
\lambda_{1,3}&=[1-m_1(16l_2^2+m_1^2+24l_2+8)^{-1/2}]/2,\\
f_{2,2}& = (2l_2^2+3l_2+1)x^2+((2l_2+1)y+1)(y-1),\\
\lambda_{2,2}&=1.\\
\end{aligned}
\end{equation*}

$\bullet$ Family $\mathcal{V}_7$ is equivalent to $\mathcal{V}_3,$ just by interchanging the left and right differential systems.

$\bullet$ For the last family $\mathcal{V}_8$ we can take
\begin{equation*}
\begin{aligned}
f_{1,1}&=f_{2,1} =1-y,\\
f_{1,2}& =-(l_2+1)l_2x^2/2+(l_2y+1)(1-y),\\
\lambda_{1,1}&=\lambda_{2,1}-1=l_2-1,\\
\lambda_{1,2}&=1,\\
f_{2,2}& = \left(-m_2+(4l_2^{2}+m_2^{2}+4l_2)^{1/2} \right) x/2+l_2y+1,\\
f_{2,3}& = \left(-m_2-(4l_2^{2}+m_2^{2}+4l_2)^{1/2}\right) x/2+l_2y+1,\\
\lambda_{2,2} &=(1+(4l_2^2+m_2^2+4l_2)^{-1/2}m_2)/2,\\
\lambda_{2,2} &=(1-(4l_2^2+m_2^2+4l_2)^{-1/2}m_2)/2,\\
\end{aligned}
\end{equation*}
being  $H_j(0,y)=(1-y)^{l_2}(l_2y+1),$ for $j=1,2.$ 
\end{proof}

\section{The Maximal Order of a Weak-Focus and the Bifurcation of Crossing Limit Cycles}\label{se:analytic}	

In this section, we will provide the conditions of the parameters such that systems \eqref{eq:1_1c}, \eqref{eq:1_1}, \eqref{eq:1_1ac} have the maximal order of a weak-focus located at the origin and also the unfolding of crossing limit cycles of small amplitude in each family. This is done in  Propositions~\ref{horc}, \ref{hor} and \ref{verc}. The complete study of system \eqref{eq:1_1a} is more intricate.  Proposition~\ref{ver} provides the maximal order of each weak-focus and some values of the parameters such that this maximality is attained. Finally, in Proposition~\ref{prop45} we get the complete unfolding of some of them. Consequently, the proof of Theorem~\ref{thmmainn} is finished. Although the proofs of Theorems~\ref{Thm:continuo4} and \ref{Thm:continuo5} will be done in the next section, the explicit unfoldings follow from the following results. This section is structured in two subsections. The first contains all the results referred to the case with the $x$-axis as the separation straight line. The second is devoted to the result being the $y$-axis as the separation straight line.

As we have explained in the introduction, in the following results we will always have one crossing limit cycle of small amplitude less than the order of each weak-focus. Because our canonical forms have no sliding segment.

\subsection{The Horizontal Case} 
\begin{pro}\label{horc} 
The maximal weak-focus order of the origin of the piecewise differential system \eqref{eq:1_1c} is $4.$ This maximal property is obtained when the parameters are on
\begin{equation*}
\mathcal{T}^c=\{d=2l+n_1+n_2= m_1-m_2=0; m_2(n_1-n_2)(n_1+n_2-1)\neq0\}.
\end{equation*}
Additionally, the weak-foci on $\mathcal{T}^c$ unfold $3$ limit cycles of small amplitude bifurcating from the origin, multiplicities taken into account and perturbing inside family \eqref{eq:1_1c}.
\end{pro}

\begin{proof}
The first necessary condition to have a nondegenerate equilibrium point of center-focus type at the origin of \eqref{eq:1_1c} is $d=0.$ Because the trace and the determinant of the Jacobian matrix are zero and one, respectively. With the mechanism described in Section~\ref{se:preliminaries}, straightforward computations provide the first Lyapunov quantities $W_n.$ In particular, $W_1=0$ because $d=0$ and 
\begin{equation}\label{csthc}
\begin{aligned}
W_2=&2(m_1-m_2)/3,\\
W_3=&\pi m_2(2l+n_1+n_2)/8,\\
W_4=&4m_2(n_1-n_2)(6l+4n_1+4n_2-1)/45.
\end{aligned}
\end{equation}
The proof of the maximality follows checking that $W_2(\mathcal{T}^c)=W_3(\mathcal{T}^c)=0$ and  $W_4(\mathcal{T}^c)=4m_2(n_1-n_2)(n_1+n_2-1)/45\ne0$ and that the solutions of the polynomial system $\{W_2=W_3=W_4=0\}$ provide the centers detailed in Corollary~\ref{coroh}, which are centers using Proposition~\ref{cent3} and the continuity condition.

As the determinant of the Jacobian matrix of $W_2,W_3$ with respect to $(m_1,n_2)$ on $\mathcal{T}^c,$
\begin{equation*}
\det J=\det(\text{Jac}\left[(W_2,W_3),(m_1,n_2)\right]|_\mathcal{T^c})= \left| \begin {array}{cc} 2/3&0\\ \noalign{\medskip}0&1/8\,\pi \,{
	\it m_2}\end {array} \right|=\pi m_2/12,
\end{equation*}
is different from zero, we have two hyperbolic limit cycles bifurcating from the origin under the condition $d=0.$ The third limit cycle emerges from the origin in a similar way as the classical Hopf bifurcation being $d$ small enough and different from zero. As we have explained previously. The unfolding taking into account the multiplicities can be proved using the results in \cite{HanYan2021}.
\end{proof}

\begin{pro}\label{hor} 
The maximal weak-focus order of the origin of the piecewise differential system \eqref{eq:1_1} is $5.$ This maximality is obtained when the parameters are on
\begin{equation*}
\mathcal{T}=\{d_1=d_2=2l_2+3n_1-n_2= 2l_1-n_1+3n_2=m_1- m_2=0; m_2(n_1+n_2-1)(n_2-n_1)\neq0\}.
\end{equation*}
Additionally, the weak-foci on $\mathcal{T}$ unfold $4$ limit cycles of small amplitude bifurcating from the origin, multiplicities taken into account, and perturbing inside family \eqref{eq:1_1}.
\end{pro}

\begin{proof} 
As the proof follows similarly to the proof of Proposition~\ref{horc}, we only detail the differences. For system \eqref{eq:1_1}, the origin is a nondegenerate weak-focus when $d_1=d_2=0.$ The first Lyapunov quantities are	
\begin{equation}\label{csth}
\begin{aligned}
W_2  =&2(m_1-m_2)/3,\\
W_3  =&\pi m_2(l_1+l_2+n_2+n_1)/8,\\
W_4=&-2m_2(l_1+l_2+1)(3l_1+l_2+4n_2)/45,\\
W_5= & -\pi m_2(l_1+l_2+1)(l_1-l_2)^2/1536, \\
W_6=& -`2m_2(l_2+2)(l_2-1)(l_1+l_2+1)(l_1-l_2)/4725.\\
\end{aligned}
\end{equation}
Straightforward computations show that over $\mathcal{T}$  we have $W_2=W_3=W_4=0$ and $W_5=\pi m_2(n_1+n_2-1)(n_1-n_2)^2/384\ne0.$ The maximality follows from the fact that  
$W_6^2\subset\langle W_1,W_2,\dots,W_5\rangle$ and that under the conditions $W_2=W_3=W_4=W_5=0$ we have the centers detailed in Theorem~\ref{centro1}, which are centers because of Proposition~\ref{cent3}.

The unfolding of limit cycles bifurcating from the origin follows also similarly to the previous proof. When $d_1=d_2=0$ the determinant of the Jacobian matrix of $W_2,W_3,W_4$ with respect to $(l_2,m_1,n_2)\ne0$ over $\mathcal{T}$ is 
$\text{det}(J)=\pi(m_2)^2(n_1+n_2-1)/90,$
where 
\begin{equation*}
J=\left(\begin{array}{ccc} 0&\dfrac{2}{3}&0\\[10pt]
\dfrac{\pi}{8}m_2&\dfrac{\pi}{8}m_2&\dfrac{\pi}{8}m_2\\[10pt]
\dfrac{2}{45}m_2( 21n_1-19n_2-1)&0&\dfrac{4}{45}m_2(9n_1-11n_2+1)\end {array} \right).
\end{equation*}
Hence, we have three limit cycles of small amplitude and the fourth bifurcates taking $d_2=0$ and $d_1\ne0$ small enough but with an adequate sign. The unfolding taking into account the multiplicities is proved using \cite{HanYan2021}.
\end{proof}

\begin{proof}[Proof of Theorem~\ref{centro1}] The necessary conditions for having a center at the origin are $d_1=d_2=0$ and the equation \eqref{csth} obtained in the proof of Proposition~\ref{hor}. The sufficiency is provided by Proposition~\ref{cent3}.
\end{proof}

\subsection{The Vertical Case}

\begin{pro}\label{verc} 
The maximal weak-focus order of the origin of the piecewise differential system \eqref{eq:1_1ac} is $3.$ This maximality is obtained when the parameters are on
\begin{equation*}
\mathcal{F}^c=\{d_1=d_2=l_1-l_2=0; (m_1+m_2)(l_2+n)\neq0\}.
\end{equation*}
Additionally, the weak-foci on $\mathcal{F}^c$ unfold $2$ limit cycles of small amplitude bifurcating from the origin, multiplicities taken into account and perturbing inside family \eqref{eq:1_1ac}.
\end{pro}

\begin{proof}
The proof follows basically using the same steps as the proof of Proposition~\ref{horc}. Here for computing th Lyapunov quantities first we need to consider a rotation of angle $-\pi/2$ in order that the separation straight line be the  $x$-axis. Once again we have that $W_1=0,$ when $d_1=d_2=0.$ The first Lyapunov quantities are
\begin{equation}\label{csthvc}
\begin{aligned}
W_2  =&4(l_1-l_2)/3,\\
W_3  =&\pi(m_1+m_2)(l_2+n)/8.\\ 
\end{aligned}
\end{equation}
When $W_2=W_3=0$ we have one center at the origin as the ones listed in Corollary~\ref{centvc}, but they are centers because of Proposition~\ref{cent4}, assuming the continuity condition. Consequently, the property of maximality and the existence of the condition $\mathcal{F}^c$ follow. Like in the previous two proofs, the complete unfolding also follows. Here the linearity condition of $W_2$ with respect to $l_1$ or $l_2$ provides the first limit cycle of small amplitude. The second, as above, taking $d_2=0$ and $d_1\ne0$ small enough.	
\end{proof}

\begin{pro}\label{ver} 
The maximal weak-focus order of the origin of the piecewise differential system \eqref{eq:1_1a} is $8.$ In particular, there are at least four families exhibiting this maximality:
\begin{equation}\label{eq:famF}
\begin{aligned}
\mathcal{F}_1=&\{d_1=d_2=m_1=m_2= 0,l_1=-13/4,l_2=-3/2,n_2=-1/2\},\\
\mathcal{F}_2=&\{d_1=d_2=m_1=m_2= 0,l_1=9/4,l_2=1/2,n_2=3/2\},\\
\mathcal{F}_3^\pm=&\{d_1=d_2=m_1=m_2= 0,\\
&\quad l_1=(\pm 1+6 l_2+f(l_2))/4,n_2= 1/2\pm f(l_2)/5,l_2\not\in \mathcal{L}\},
\end{aligned}
\end{equation}
where $\mathcal{L}=\{-3/4, -1/4,-3/2\} $ and $f(l_2)=\sqrt{20l_2^2+20l_2+10}.$
\end{pro}

\begin{proof} 
The proof follows basically using the same steps as the previous proofs, but the computations are more intricate. As above, we will start assuming $d_1=d_2=0$ to get $W_1=0.$ Next, in order to apply the algorithm described in Section~\ref{se:preliminaries}, as in the proof of Proposition~\ref{verc}, we need to do a rotation of angle $-\pi/2$ to compute the Lyapunov quantities in this case. As usual, the property of maximality will follow solving the algebraic system of equations
\begin{equation}\label{s7}
\mathcal{S}_7=\{W_2=W_3=W_4=W_5=W_6=W_7=0\},
\end{equation}
checking that there exists at least one real solution such that $W_8\neq0$ and proving that all the solutions of
\begin{equation}\label{s8}
\mathcal{S}_8=\{W_2=W_3=W_4=W_5=W_6=W_7=W_8=0\},
\end{equation}
imply $W_k=0=0$ for $k\geq9.$ This last step is a consequence of Proposition~\ref{cent4}. Finally, we will prove the unfolding of $7$ limit cycles described in the last statement, using, in this last step, the parameters $d_1,d_2.$

Straightforward computations allow us to get the first Lyapunov quantities which are polynomials in the parameters space $(l_1,l_2,m_1,m_2,n_1,n_2).$ Because of the size, we only detail the first one which, using $W_2=0,$ provides the condition
\begin{equation}\label{n1}
n_1=-2l_1+2l_2+n_2.
\end{equation} 
The direct application of the algorithm of Section~\ref{se:preliminaries} provides the coefficients of the displacement function \eqref{eq:11} that write, some of them, as polynomials in $\pi,$ before using that the previous should vanish. So we have
\begin{equation}\label{eq:8liap}
\begin{aligned}
W_3=&\pi\widetilde{W}_3, \quad W_4=\widetilde{W}_4, \quad W_5=\pi\widetilde{W}_5, \quad W_6=\widetilde{W}_6, \quad W_7=\widetilde{W}^{[0]}_7+\pi\widetilde{W}^{[1]}_7,\\
W_8=&\widetilde{W}^{[0]}_8+\pi\widetilde{W}^{[1]}_8, \quad W_9=\widetilde{W}^{[0]}_9+\pi\widetilde{W}^{[1]}_9+\pi^2\widetilde{W}^{[2]}_9, \\ W_{10}=&\widetilde{W}^{[0]}_{10}+\pi\widetilde{W}^{[1]}_{10}+\pi^2\widetilde{W}^{[2]}_{10}, \quad W_{11}=\widetilde{W}^{[0]}_{11}+\pi\widetilde{W}^{[1]}_{11}+\pi^2\widetilde{W}^{[2]}_{11}+\pi^3\widetilde{W}^{[3]}_{11},
\end{aligned}
\end{equation}
being $\widetilde{W}_i^j$ polynomials with rational coefficients in $(l_1,l_2,m_1,m_2,n_2).$ Using a computer algebra system  we can see that $\widetilde{W}^{[0]}_7\in\langle \widetilde{W}_3\dots,\widetilde{W}_{6}\rangle$, $\widetilde{W}^{[1]}_8\in\langle \widetilde{W}_3,\dots,\widetilde{W}_{6},\widetilde{W}^{[1]}_7\rangle,$ and $\widetilde{W}_9^{[i]},  \big( \widetilde{W}^{[0]}_{10}\big)^2, \widetilde{W}^{[1]}_{10},$ $\widetilde{W}^{[2]}_{10}, \widetilde{W}_{11}^{[i]}\in \langle \widetilde{W}_3\dots,\widetilde{W}_{6},\widetilde{W}^{[1]}_7,\widetilde{W}^{[0]}_8\rangle.$ 
Moreover, we can write
\begin{equation*}
\begin{aligned}
\widetilde{W}_3=&(l_2m_2-l_1m_1+m_2n_2+m_1n_2+2l_2m_1)/8,\\
\widetilde{W}_4=&(96l_2^3-240l_2^2l_1+144l_2^2n_2+192l_2l_1^2-240l_2l_1n_2-8l_2m_2^2\\
&-8l_2m_2m_1+60l_2n_2^2-48l_1^3+96l_1^2n_2-60l_1n_2^2-8m_2^2n_2\\
&-8m_2m_1n_2-48l_2^2+72l_2l_1-36l_2n_2-24l_1^2+36l_1n_2)/45,\\
\widetilde{W}_5=&(l_1-l_2)(22l_2^2m_2+52l_2^2m_1+67l_2l_1m_2-52l_2l_1m_1+45l_2m_2n_2\\
&-52l_1^2m_2+13l_1^2m_1+15l_1m_2n_2+10m_2n_2^2+33l_2m_2+26l_2m_1\\
&-26l_1m_2-13l_1m_1+20m_2n_2)/384,\\
\widetilde{W}_6=&8(l_2-l_1)(15664l_2^4-41576l_2^3l_1+41128l_2^3n_2+40440l_2^2l_1^2\\
&-79808l_2^2l_1n_2+41730l_2^2n_2^2-17272l_2l_1^3+50992l_2l_1^2n_2-51084l_2l_1n_2^2\\
&+20038l_2n_2^3+2744l_1^4-10776l_1^3n_2+15582l_1^2n_2^2-11322l_1n_2^3\\
&+3920n_2^4-10024l_2^3+17980l_2^2l_1-19942l_2^2n_2-10964l_2l_1^2\\
&+22268l_2l_1n_2-13869l_2n_2^2+2240l_1^3-6250l_1^2n_2+7023l_1n_2^2-3482n_2^3\\
&+1096l_2^2-1416l_2l_1+1219l_2n_2+434l_1^2-557l_1n_2+444n_2^2)/20475,\\
\widetilde{W}^{[1]}_7=&m_1(l_2-l_1)\left[  n_2l_2(2l_2+1)(l_2+1)\widetilde{W}_{7,20}(l_2)+(2l_2-l_1)\widetilde{W}_{7,23}(l_2)\right],\\
\widetilde{W}^{[0]}_8=& (l_2-l_1)\left[  n_2\widetilde{W}_{8,14}(l_2)+(2l_2-l_1)(2l_2-2l_1-1)\widetilde{W}_{8,24}(l_2)\right],\\
\widetilde{W}^{[0]}_{10}=&\frac{1024}{21049875}\prod_{i=1}^{9}\widetilde{\mathcal{R}}_i,  
\end{aligned}
\end{equation*}
where 
\begin{equation*}
\begin{aligned}
\widetilde{\mathcal{R}}_1=& l_2, \quad \widetilde{\mathcal{R}}_2=l_2-2l_1-1, \quad \widetilde{\mathcal{R}}_3=2l_2-l_1+1, \quad \widetilde{\mathcal{R}}_4=l_2-l_1,\quad \widetilde{\mathcal{R}}_5=l_2+1,\\
\widetilde{\mathcal{R}}_6=& 2l_2-l_1, \quad \widetilde{\mathcal{R}}_7=2l_2-2l_1+1, \quad \widetilde{\mathcal{R}}_8=2l_2-2l_1-1, \quad \widetilde{\mathcal{R}}_9=l_2-2l_1,
\end{aligned}
\end{equation*}
and $\widetilde{W}_{7,k}$ and $\widetilde{W}_{8,k}$ are polynomials with rational coefficients of degree $k$. From the above computations, systems \eqref{s7} and \eqref{s8}, using the condition \eqref{n1}, are now written as
\begin{equation*}
\begin{aligned}
\widetilde{\mathcal{S}}_7&=\{\widetilde{W}_3=\widetilde{W}_4=\widetilde{W}_5=\widetilde{W}_6=\widetilde{W}^{[1]}_7=0\},\\
\widetilde{\mathcal{S}}_8&=\{\widetilde{W}_3=\widetilde{W}_4=\widetilde{W}_5=\widetilde{W}_6=\widetilde{W}^{[1]}_7=\widetilde{W}^{[0]}_8=0\}.
\end{aligned}
\end{equation*}
As $\big( \widetilde{W}^{[0]}_{10}\big)^2\in \langle \widetilde{W}_3\dots,\widetilde{W}_{6},\widetilde{W}^{[1]}_7,\widetilde{W}^{[0]}_8\rangle$ it is clear that 
$\big( \widetilde{W}^{[0]}_{10}\big)^2=0$ on $\widetilde{\mathcal{S}}_8$ and so also $\widetilde{W}^{[0]}_{10}=0.$ But although it is not necessary to use $\widetilde{W}^{[0]}_{10}$ to solve the center problem, it is useful to use it.  Hence, after considering the equivalent system
\begin{equation*}
\widetilde{\mathcal{S}}_{10}=\{\widetilde{W}_3=\widetilde{W}_4=\widetilde{W}_5=\widetilde{W}_6=\widetilde{W}^{[1]}_7=\widetilde{W}^{[0]}_8=\widetilde{W}^{[0]}_{10}=0\},
\end{equation*}
we obtain the families of the statement of Theorem~\ref{centro2}. So, with Proposition~\ref{cent4} we have that all are center families and, consequently, the maximal weak-focus order is $8$ and the first statement follows.

The second part of the statement follows solving partially system $\widetilde{\mathcal{S}}_7.$ More concretely, solving
\begin{equation*}
\widetilde{\mathcal{S}}_7^{[1]} =\{\widetilde{W}_3=\widetilde{W}_4=\widetilde{W}_5=\widetilde{W}_6=m_1=0\}
\end{equation*}
and obtaining weak-foci of order $8$. In fact, we have that on each of them $W_i=0,$ for $i=3,\dots,7,$ and $W_8\ne0.$ More concretely, 
\begin{equation*}
\begin{aligned}
W_8(\mathcal{F}_1)&= 2/3,\\
W_8(\mathcal{F}_2)&=-2/3, \\
W_8(\mathcal{F}_3^\pm)&=-\frac{(4l_2+3)^2(4l_2+1)^2}{189000}\left[
232l_2^3+ 348l_2^2+ 222l_2 + 53 \pm(52l_2^2+52l_2+17)f(l_2)\right].
\end{aligned}
\end{equation*}
\end{proof}

\begin{pro}\label{prop45}
The weak-foci $\mathcal{F}_1$ and $\mathcal{F}_2$ defined in \eqref{eq:famF} unfold $7$ limit cycles of small amplitude bifurcating from the origin, multiplicities taken into account, when we perturb inside family \eqref{eq:1_1a}.
\end{pro}

\begin{proof}
We will follow the same unfolding procedure as in the previous results assuming $d_2=0.$ We will focus our attention only to the point $\mathcal{F}_1,$ the other follows similarly. Using the linearity dependence on $d_1$ and $n_1$ of $W_1$ and $W_2$ defined in \eqref{n1}, we can restrict our analysis to the study of the transversality condition of the Taylor series of the next Lyapunov quantities near $\mathcal{F}_1,$ with respect to the parameters $(m_1,m_2,l_1,l_2,n_2).$ Taking the perturbation
\begin{equation*}
\mathcal{F}_{1,\e}=\left[m_1=\e_1, m_2=\e_2, l_1=-13/4+\e_3,l_2=-3/2+\e_4,n_2=-1/2+\e_5\right]
\end{equation*}
and with the linear change of variables in the parameter space,
\begin{equation*}
\begin{aligned}
\e_1&=\dfrac{5}{4}u_7+\frac{5}{8}u_6+\frac{3215}{672}u_4,\quad \e_2=\dfrac{25}{8}u_7+\dfrac{25}{16}u_6+\frac{14635}{1344}u_4,\\
\e_3&=\dfrac{14}{3}u_1+\frac{64}{35}u_3,\quad \e_4= -\frac{208}{3}u_3-\frac{512}{35}u_5, \quad \e_5 = u_7,
\end{aligned}
\end{equation*}
we have that the Taylor series of the Lyapunov quantities write as 
\begin{equation*}
\begin{aligned}
W_j(u)=& u_j+O(u^2), \text{ for } j=3,\ldots, 6, \\
W_7(u)=&\frac{5470997}{2496}\pi u_1-\frac{184}{9} u_2+\frac{781571}{1690}\pi u_3-12u_4 +O(u^{2}),\\
W_8(u)=& 2/3+O(u).
\end{aligned} 
\end{equation*} 
We notice that if $u_i=0$ we have a weak-focus of order $8$. Moreover, with the Implicit Function Theorem we have new coordinates $v_3,\ldots,v_7$, in the parameter space, such that $W_j(u)=v_j,$ para $j=3,\ldots,6$ and $u_7=v_7.$ Hence  the transversality condition is satisfied up to $W_6.$ The last step is the computation of the Taylor series of $W_7$ when $v_3=v_4=v_5=v_6=0.$ Straightforward computations provide \[
W_7(v_7)=-\frac{80339}{31104} v_7^3+O(v_7^4).
\]
The unfolding is complete because the above first coefficient has an odd power in the remaining parameter $v_7.$ More details on the used technique can be seen in \cite{FerGinTor2021}. 
\end{proof}

We remark that the complete unfolding study for the other families of weak-foci in Proposition~\ref{ver} is more difficult because of the dependence on the parameter $l_2.$ But it can be seen that only linear developments are not enough.
   
\begin{proof}[Proof of Theorem~\ref{centro2}] This proof is a direct consequence of the proofs of Propositions \ref{ver} and \ref{cent4}, since every candidate to be a center is given nullifying the first eight Lyapunov quantities $W_i$, for $i=1,\ldots,8$ obtained in \eqref{eq:8liap}.
\end{proof}

\section{The cyclicity problem in the continuous classes}
We finish the work by studying the maximum number of limit cycles of small amplitude that bifurcate from the origin in the continuity classes \eqref{eq:1_1c} and \eqref{eq:1_1ac}. That is obtaining its cyclicity and proving Theorems~\ref{Thm:continuo4} and \ref{Thm:continuo5}. 

\begin{proof}[Proof of Theorem~\ref{Thm:continuo4}] Using Proposition~\ref{horc} we know that from the origin of system \eqref{eq:1_1c} bifurcates at least $3$ limit cycles of small amplitude. The multiplicity property follows like the previous results using \cite{HanYan2021}. The upper bound follows from Theorem~9, of Chapter~2, given in \cite{Roussarie1998} because the ideal, $I=\langle W_2,W_3,W_4\rangle,$ generated by the Lyapunov quantities, given in \eqref{csthc}, is radical. The radicality proves that under the condition $d=0$ we have at most two limit cycles. The third limit cycle appears, using $d$, as in a classical Hopf bifurcation. See more details in \cite{AndLeoGorMai1973} or again \cite{Roussarie1998}.
\end{proof}

\begin{proof}[Proof of Theorem~\ref{Thm:continuo5}] The proof follows analogously as the above proof using that the ideal generated by the Lyapunov quantities $W_2,W_3$ defined in \eqref{csthvc} is also radical.
\end{proof}

The above approach can not be used for studying the cyclicity of \eqref{eq:1_1} and \eqref{eq:1_1a} because the ideal generated by the corresponding Lyapunov quantities it is not radical.

\section*{Acknowledgements}

This work has been realized thanks to the Brazilian CAPES Agency (Coordena\c{c}\~ao de Aperfei\c{c}oamento de Pessoal de N\'{\i}vel Superior - Finance Code 001), the Catalan AGAUR Agency (grant 2017 SGR 1617), the Spanish Ministerio de Ci\'encia, Innovaci\'on y Universidades via the Agencia Nacional de Investigaci\'on (grants PID2019-104658GB-I00 and CEX2020-001084-M), and the European Union's Horizon 2020 research and innovation programme (grant Dynamics-H2020-MSCA-RISE-2017-777911).

\bibliographystyle{abbrv}
\bibliography{biblio} 
\end{document}